# ΜΑΘΗΜΑΤΙΚΗ 89
# ΕΠΙΘΕΩΡΗΣΗ

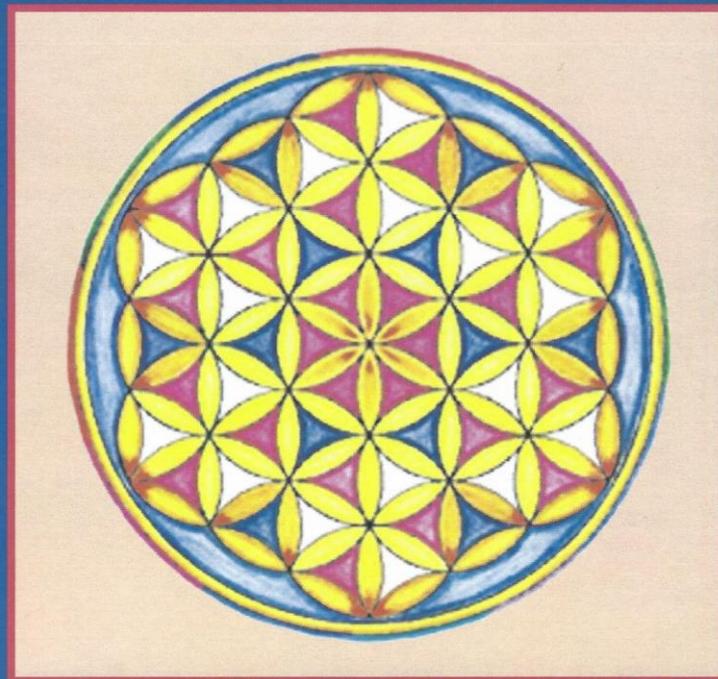

Ελληνική Μαθηματική Εταιρεία

**To cite this article:**

Ρίζος, Ι. (2018). Φιλοσοφικές και επιστημολογικές θεωρήσεις της έννοιας του χώρου. *Μαθηματική Επιθεώρηση*, 89, 43-60.

# Φιλοσοφικές και επιστημολογικές θεωρήσεις της έννοιας του χώρου

**Ιωάννης Ρίζος**
Δρ Μαθηματικών Πανεπιστημίου Πατρών
irizos@hotmail.gr

**Περίληψη**
Από την αρχαιότητα η εννοιολογική αντίληψη του χώρου μεταβλήθηκε επίπονα και με σχετικά αργούς ρυθμούς. Πέρασε μέσα από μυθολογικές περιγραφές, θρησκευτικές δοξασίες, μεταφυσικές κοσμοθεωρίες και κοσμολογικά πρότυπα με μηχανιστική δομή, μέχρι να φτάσει στη νευτώνεια αντίληψη του απέραντου και ομογενούς τρισδιάστατου συνεχούς, που κυριαρχεί μέχρι τις ημέρες μας ως "κοινή αντίληψη" ακόμα και στην εκπαίδευση. Αυτή η αντίληψη για τον χώρο η οποία σήμερα φαίνεται αυτονόητη, αποτελεί την κατάληξη μιας μακράς ιστορικής εξέλιξης, την οποία θα προσπαθήσουμε να σκιαγραφήσουμε. Στόχος μας είναι να παρουσιάσουμε συνοπτικά την εννοιολογική πορεία των προτύπων σχετικά με τον φυσικό χώρο από την αρχαιότητα μέχρι τον 20ο αιώνα, εστιάζοντας στις θεωρίες και τα μοντέλα που επιδέχονται μαθηματικοποίηση. Ταυτόχρονα επιχειρούμε να αναδείξουμε μια από τις σημαντικότερες αφετηρίες των γεωμετρικών εννοιών, που είναι η ανάγκη μαθηματικής μοντελοποίησης του φυσικού χώρου, καθώς στη σύγχρονη εκπαίδευση απουσιάζει μια τέτοια θεώρηση.

### 1. Πρώιμες αναπαραστάσεις του κόσμου

Για πολλούς αιώνες οι αναπαραστάσεις της αρχαίας Αστρονομίας καθόριζαν και την αντίληψη του χώρου. Ολόκληρος ο κόσμος περιστρεφόταν φυσικά και αβίαστα γύρω από τον ίδιο σταθερό άξονα ή το ίδιο σταθερό σημείο. Οι πρώτες κοσμολογικές αντιλήψεις δεν ήταν απαλλαγμένες από μυθολογικά, ανιμιστικά και ανθρωπομορφικά στοιχεία, που παρατηρούνται ακόμα και στη σκέψη των παιδιών (Πατρώνης & Ρίζος 2011, σ. 550).

Για παράδειγμα, η κοσμολογία των Πυθαγορείων υπήρξε αδιαχώριστο κράμα ορθολογισμού και μυθοπλασίας, και ήταν αλληλένδετη με τα Μαθηματικά, τη Μουσική και τον μυστικισμό. Τα ουράνια σώματα δεν ήταν απλά άψυχα στολίδια καρφωμένα στον ουράνιο θόλο. Ήταν θεία, αιώνια όντα, είχαν ψυχή και νόηση, καθώς και πλήρη αίσθηση της μουσικής αρμονίας. Έτσι επέλεξαν να υιοθετήσουν το τελειότερο σχήμα –τη σφαίρα– και να περιφέρονται αιώνια γύρω από τη Γη ακολουθώντας την τελειότερη τροχιά –τον κύκλο– σε αποστάσεις αντίστοιχες με τα αρμονικά διαστήματα της μουσικής κλίμακας (Van der Waerden 1961, 2003). Οι πέντε πλανήτες που από την αρχαιότητα παρατηρούνταν με γυμνό μάτι (Ερμής, Αφροδίτη, Άρης, Δίας, Κρόνος) μαζί με τη Γη, τη Σελήνη και τον Ήλιο έφταναν τον αριθμό οκτώ. Προκειμένου, κατά πάσα πιθανότητα, να συμπληρωθεί η ιερή δεκάδα, ο Πυθαγόρειος Φιλόλαος[1] επινόησε μια αόρατη "αντι-Γη" με μάζα ίση προς αυτή της Γης, την *Ἀντίχθωνα*, καθώς και ένα αιώνιο κεντρικό πυρ, την *Ἑστία*, γύρω από την οποία περιφερόταν ολόκληρος ο κόσμος, μηδέ της Γης εξαιρουμένης (Αριστοτέλης, *Περί Ουρανού*, Β 293a 18; Αέτιος ΙΙ, 7, 7 [DK 44a 16]). Είναι βεβαίως αξιοσημείωτη η λειτουργία αυτής της πρώιμης μαθηματικοποίησης: το μοντέλο δεν δημιουργείται με αποκλειστικό σκοπό να ερμηνεύσει μια παρατηρούμενη πραγματικότητα, αλλά αντίθετα, η πραγματικότητα εξωθείται να ικανοποιήσει τις μυστικιστικές απαιτήσεις του κοσμοθεωρητικού προτύπου!

Σύμφωνα με τον Αριστοτέλη (*Περί Ουρανού*, Β 293a 23), το μοντέλο του Φιλολάου ήταν σε θέση να εξηγήσει το φαινόμενο της εναλλαγής ημέρας-νύχτας, και ενδεχομένως ο δημιουργός του να ήθελε επιπλέον να ερμηνεύσει την προέλευση του ηλιακού φωτός. Ακόμη, όπως διασώζει ο Αριστοτέλης, οι Πυθαγόρειοι ήσαν οι πρώτοι που διατύπωσαν την άποψη ότι στο κέντρο του σύμπαντος υπάρχει φωτιά[2] και ότι είναι λάθος να τοποθετείται η γη στο κέντρο του σύμπαντος,[3] θεωρώντας πως η πιο τιμητική θέση αρμόζει στο πιο πολύτιμο πράγμα: στη φωτιά που είναι πολυτιμότερη από τη γη[4]. Με αυτή την «Ἑστίαν τοῦ παντός» και του «Διὸς οἶκον» που συλλαμβάνει ο Φιλόλαος (Αέτιος ΙΙ, 7, 7 [DK 44a 16]) προαναγγέλλει την μετέπειτα ηλιοκεντρική θεωρία του Αρίσταρχου. Θα πρέπει ωστόσο να ληφθεί υπόψη πως ένα μάλλον βαθύτερο θρησκευτικό αίσθημα οδηγεί τον Φιλόλαο να διατυπώσει τη θεώρησή του. Φυσικά η θεώρηση αυτή δεν παύει να αποτελεί ένα πρώτο βήμα αμφισβήτησης του γεωκεντρισμού, καθώς στηρίζεται περισσότερο στη *νόηση* και παρεκκλίνει από την αμεσότητα της αίσθησης για την (φαινομενική) περιστροφή της ουράνιας σφαίρας.

---

[1] Ο Van der Waerden (1961, σ. 7) αποδίδει την ιδέα αυτή στον Ικέτα τον Συρακούσιο, όπως ακριβώς και ο Κοπέρνικος στο χειρόγραφο του έργου του *De revolutionibus orbium coelestium*, καθώς και στις μεταγενέστερες εκδόσεις αυτού.
[2] «ἐπὶ μὲν γὰρ τοῦ μέσου πῦρ εἶναί φασι», Αριστοτέλης, *Περί Ουρανού*, Β 293a 21.
[3] «μὴ δεῖν τῇ γῇ τὴν τοῦ μέσου χώραν ἀποδιδόναι», ό.π., Β 293a 28.
[4] «Τῷ γὰρ τιμιωτάτῳ οἴονται προσήκειν τὴν τιμιωτάτην ὑπάρχειν χώραν, εἶναι δὲ πῦρ μὲν γῆς τιμιώτερον», ό.π., Β 293a 30.

Οι παραπάνω απόψεις και δοξασίες είχαν ως συνέπεια να επικρατήσει η αντίληψη ότι το σύμπαν είναι κλειστό και πεπερασμένο. Ο *χώρος* περιορίστηκε στα μεσοδιαστήματα μιας σειράς περιστρεφόμενων ομόκεντρων σφαιρών, στην εξώτατη από τις οποίες βρίσκονταν οι απλανείς αστέρες και στις υπόλοιπες οι πλανήτες, ο Ήλιος και η Σελήνη. Προφανώς γι' αυτό τον λόγο το κυριότερο σχήμα στην αρχαία ελληνική γεωμετρία ήταν ο κύκλος. Αντίθετα, η δυνατότητα ενός άπειρου σύμπαντος συναντάται στην κοσμολογία των Ιώνων φιλοσόφων και ιδιαίτερα αυτών της σχολής της Μιλήτου, καθώς και στις αντιλήψεις του Δημόκριτου και του Επίκουρου.

### 2. Η πλατωνική *Χώρα* και ο αριστοτελικός *Τόπος*

Σήμερα, ως απαρχή της επιστημονικής Αστρονομίας θεωρείται το γεωκεντρικό μοντέλο του Εύδοξου του Κνίδιου και του Κάλλιππου του Κυζικηνού, το οποίο χαρακτηρίζεται από τον συνδυασμό της θεωρίας με την παρατήρηση (Dreyer 1953, σ. 107). Η σχέση του Εύδοξου και του Κάλλιππου με τον Πλάτωνα και τον Αριστοτέλη είναι γνωστή,[5] επομένως δεν προκαλεί εντύπωση ότι τον 4ο αιώνα π.Χ., με επίκεντρο την Ακαδημία, ξεκινά μια συντονισμένη προσπάθεια γεωμετρικοποίησης *του χώρου που μας περιβάλλει*, απαλλαγμένη από μυθολογικές δοξασίες όπως π.χ. η περιδίνηση των διάπυρων λίθων μέσα στον αιθέρα κατά τον Αναξαγόρα (Ιππόλυτος, *Έλεγχος* Ι, 8, 3-10 [DK 59a 42]).[6] Το πρόγραμμα αυτό, βασισμένο σε μια θεωρία κίνησης επί ομαλώς περιστρεφόμενων σφαιρών, έμελλε να σημαδέψει την κοσμολογική σκέψη για πολλούς αιώνες (Kuhn 1957, σ. 59).

---

[5] Πρόκλος, *Σχόλια εις Ευκλείδην*, 67,2 κ.ε. και Σιμπλίκιος, *Σχόλια εις Περί Ουρανού του Αριστοτέλους*, 7,32,16 κ.ε. Για το μοντέλο των ομόκεντρων σφαιρών των Ευδόξου–Κάλλιππου βλ. Yavetz, I. (1998). On the Homocentric Spheres of Eudoxus. *Archive for History of Exact Sciences*, 52, 221-278 και Αριστοτέλης, *Μετά τα Φυσικά Λ*, 1073b 32-38.

[6] Η προ Ευδόξου (πυθαγόρειας προέλευσης) αντίληψη του Πλάτωνα ότι η φυσική πραγματικότητα για να είναι ορθολογική, οφείλει να ενσαρκώνει μια μαθηματική δομή, μεταφράζεται στον *Τίμαιο* με τις κυκλικές κινήσεις των πλανητών (38c), με την τριγωνική μορφή της ύλης (53c) και με τα πέντε κανονικά στερεά (54a). Βλ. σχετ. Κάλφας (2003).

Ειδικά για το θέμα του χώρου, ο Αριστοτέλης (*Φυσικά*, 209b17) αναφέρει ότι ο Πλάτωνας ήταν ο μόνος, μέχρι τότε, ο οποίος επιχείρησε να πει «τι είναι ο χώρος». Πραγματικά, στον *Τίμαιο* ο Πλάτωνας εκτός από τα δύο γνωστά γένη –τις *Ιδέες* και τα *Αισθητά*, που κατοικούν σε διαφορετικούς κόσμους– *κατ' ανάγκην* (Κάλφας 2003, σ. 427) εισάγει και ένα τρίτο γένος, τη *Χώρα* ή "Υποδοχή".[7] Σε πρώτη ανάγνωση δεν μπορεί να αποκλειστεί μια πρόθεση του Πλάτωνα να συνδέσει τον Κόσμο των Ιδεών με τον Φυσικό Κόσμο, κάτι το οποίο επιτυγχάνεται με την εισαγωγή της έννοιας της Χώρας ως λογική αναγκαιότητα. Ενδεχομένως, σε ένα ανώτερο τελεολογικό επίπεδο, η Χώρα να είναι απαραίτητη για να διαχωριστούν πλήρως οι Ιδέες από τα Αισθητά και να αποτραπεί η άμεση επικοινωνία τους, να προβληθούν οι μεν επάνω στα δε (χωρίς να μπορεί να συμβεί το αντίστροφο), και έτσι να "σωθούν τα φαινόμενα" από επιχειρήματα όπως αυτό του "τρίτου ανθρώπου" που εμφανίζεται στον *Παρμενίδη* (132a).

Η Χώρα υποδέχεται τον Κόσμο των Ιδεών και τον Φυσικό Κόσμο ως ένα σύνολο που δίνει "έδρα" στα στοιχεία του.[8] Δεν είναι όμως ούτε το ένα ούτε το άλλο, αφού «δεν ανήκει ούτε στην τάξη του *είδους*, ούτε στην τάξη των μιμημάτων, των εικόνων του *είδους* που έρχονται να εντυπωθούν σε αυτήν», αλλά ταυτόχρονα είναι και το ένα και το άλλο, ενσαρκώνοντας τη *μέθεξη* της πλατωνικής φιλοσοφίας (Derrida 2000). Η Χώρα φαίνεται πως είναι έτσι κάτι το "ιδεατό" και ταυτόχρονα "ζωντανό" ως *υποδοχή, μητέρα* και *τροφός*, μια ενδιάμεση οντότητα και ένα μέσο, θα λέγαμε, το οποίο περιέχει τον Κόσμο των Ιδεών, όπως και τον Φυσικό Κόσμο που φαντάζει στις αισθήσεις μας ως είδωλο του Κόσμου των Ιδεών. Η γνώση που προέρχεται από τον Φυσικό Κόσμο είναι παροδική (λόγω της αναξιοπιστίας των ανθρώπινων αισθήσεων) και μεταβλητή (λόγω των αέναων αλλαγών που υφίσταται ο φυσικός κόσμος), σε τέτοιο βαθμό που τα Αισθητά να θεωρούνται «φευγαλέες ψευδο-οντότητες» που δεν τους αξίζει ούτε να αποκαλούνται με κάποιο συγκεκριμένο όνομα (*Τίμαιος*, 49d κ.ε.). Αντίθετα, ο Κόσμος των Ιδεών, ο οποίος γίνεται αντιληπτός με τον νου, είναι *άχρονος και αναλλοίωτος*, και τα Μαθηματικά (όπως και η μαθηματική Αστρονομία) είναι ακριβώς η γνώση του κόσμου αυτού (*Πολιτεία*, 527b). Έτσι με τη διαμεσολάβηση της Χώρας, η εκάστοτε Ιδέα λειτουργεί ως το αμετάβλητο και αιώνιο πρότυπο, το αρχετυπικό μοντέλο μιας κλάσης αντικειμένων του Φυσικού Κόσμου.

---

[7] «τρίτον δὲ αὖ γένος ὂν τὸ τῆς χώρας ἀεί» Πλάτων, *Τίμαιος*, 52a8.
[8] «ἕδραν δὲ παρέχον ὅσα ἔχει γένεσιν πᾶσιν» Πλάτων, *Τίμαιος*, 52b1.

Τα παραπάνω μας οδηγούν στη διάταξη των πλατωνικών γενών με έναν τρόπο που λαμβάνει υπόψη του τη *Χώρα* (Σχήμα 2), σε αντιπαραβολή με την ευρέως διαδεδομένη πλατωνική αντίληψη περί "δυϊσμού" (Σχήμα 1).

Έτσι, η *Χώρα* καθίσταται το πλαίσιο και συνάμα το υπόβαθρο του συνόλου των φιλοσοφικών συνθηκών που θέτει ο Πλάτωνας. Το σύνολο αυτό περιλαμβάνει κυρίως τη σύνθεση του κόσμου από *τέσσερα στοιχεία* (γη, νερό, αέρας, φωτιά)[9] και την *τάξη* που επικρατεί σ' αυτόν εξαιτίας των *ίσων λόγων* μεταξύ των στοιχείων (*Τίμαιος*, 31b4-32c4). Δεδομένου λοιπόν ότι οι συνθήκες αυτές από κοινού με την πλατωνική παράδοση[10] οριοθετούσαν για αιώνες το επιστημονικό πεδίο αναζήτησης όσον αφορά την κοσμολογία, δεν θα ήταν ολότελα εσφαλμένο να θεωρήσουμε πως παίζουν τον ρόλο του *Πλατωνικού Παραδείγματος*.[11]

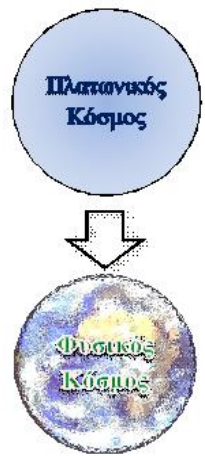
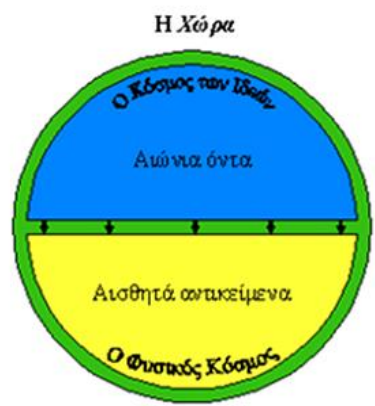

Σχήμα 1                                      Σχήμα 2

Επιχειρώντας τώρα να ορίσουμε το *Αριστοτελικό Παράδειγμα* κατ' αναλογία με το πλατωνικό, παρατηρούμε ότι οι δύο πρώτες παράμετροι, δηλαδή τα *τέσσερα στοιχεία* και η *τάξη* παραμένουν ίδιες, ενώ προστίθεται μία πέμπτη ουσία, η *πεμπτουσία*. Η σημαντική διαφορά, όμως, έγκειται στην αλλαγή του νοήματος της τάξεως. Η αριστοτελική *τάξις* δεν είναι άλλη από τον *Τόπο* ο οποίος αντιστοιχεί στην πλατωνική *Χώρα*.[12]

---

[9] Αυτά τα τέσσερα στοιχεία είναι γεωμετρικά στερεά αποτελούμενα ουσιαστικά από συνδυασμούς τριγώνων. Άρα ο φυσικός κόσμος αποτελείται από γεωμετρικά σχήματα.

[10] «ὅτι ὁ Πλάτων ταῖς οὐρανίαις κινήσεσι τὸ ἐγκύκλιον καὶ ὁμαλὲς καὶ τεταγμένον ἀνενδοιάστως ἀποδιδοὺς πρόβλημα τοῖς μαθηματικοῖς προὔτεινε, τίνων ὑποτεθέντων δι' ὁμαλῶν καὶ ἐγκυκλίων καὶ τεταγμένων κινήσεων δυνήσεται διασωθῆναι τὰ περὶ τοὺς πλανωμένους φαινόμενα», Σιμπλίκιος, *Σχόλια εις Περί Ουρανού*, 7,492,31-493,4.

[11] Στη λέξη *Παράδειγμα* προσδίδουμε την έννοια που έδωσε ο Kuhn (1962).

Οι κοσμολογικές αντιλήψεις του Αριστοτέλη ήταν άρρηκτα συνδεδεμένες με τη θεωρία του για την κίνηση των σωμάτων. Πράγματι, σύμφωνα με τη φυσική του φιλοσοφία, η Γη είναι το κέντρο του αιώνιου, κλειστού και σφαιρικού κόσμου και αποτελείται από τέσσερα στοιχεία, όπως και όλη η υποσελήνια περιοχή. Τα σώματα κινούνται μέσα σ' αυτή την περιοχή είτε, διότι μια εξωτερική δύναμη τα εξαναγκάζει,[13] είτε διότι αναζητούν ελεύθερα τη φυσική τους θέση. Πέρα από τη Σελήνη, όμως, ισχύουν εντελώς διαφορετικοί νόμοι. Εκεί, τα ουράνια σώματα, σχηματισμένα από μια quinta essentia, κινούνται ομαλώς και αενάως επάνω σε τέλειες κυκλικές τροχιές, ακριβώς εξαιτίας της τέλειας φύσης τους. Κατά συνέπεια αυτό που απομένει για να υπάρχει *τάξις* στον κόσμο, είναι να εξεταστεί η φυσική θέση, ο *τόπος* των σωμάτων στην υποσελήνια περιοχή. Με σύγχρονη ορολογία θα λέγαμε ότι ο *τόπος* λειτουργεί έτσι ως ένα είδος πλαισίου ταξινόμησης (classification framework) της ύλης. Σημειωτέον ότι ο αριστοτελικός κόσμος είναι ένθεος (πώς αλλιώς να ερμηνεύσει κανείς την ανάγκη για την ύπαρξη ενός *πρώτου κινούντος*;) και τα πάντα μέσα του διέπονται από ένα σκοπό, προς την εκπλήρωση του οποίου τείνουν.

Σχετικά με το θέμα του *τόπου* ο Αριστοτέλης ξεκινά από μηδενική βάση, προσπαθώντας να απαντήσει στο τριπλό ερώτημα αν υπάρχει ή όχι ο *τόπος*, πώς υπάρχει και ποια είναι η ουσία του (*Φυσικά*, 208a). Η άποψή του λοιπόν είναι ότι ο *τόπος* υπάρχει και μάλιστα συνδέεται με το *ὄν* (υπό την έννοια ότι κάθε *ὄν* έχει τον *τόπο* του), αλλά και με την *κίνηση* των όντων. Σημαντικό επίσης είναι το ότι ο *τόπος* «ἔχει τινὰ δύναμιν» (*Φυσικά*, 208b10) αφού κάθε σώμα αν δεν εμποδιστεί κατευθύνεται προς τον *τόπο* του.[14]

Θα πρέπει εδώ να επισημάνουμε ότι ο *τόπος* διακρίνεται από τον Αριστοτέλη σε *κοινόν* και *ἴδιον* (*Φυσικά*, 209a 32). Από τη χρήση, όμως, της λέξης *τόπος* (χωρίς τον προσδιορισμό *κοινός* ή *ἴδιος*) φαίνεται ότι αυτή δεν εκφράζει έναν γενικό χώρο, αλλά αποτελεί ένα περιορισμένο σε έκταση τμήμα του, αφού όταν γίνεται αναφορά στον κοινό *τόπο* κάποιων σωμάτων ή σε μεγάλες περιοχές (λ.χ. το χάος του Ησιόδου), χρησιμοποιείται η λέξη *Χώρα* (*Φυσικά*, 208b7 και b32 αντίστοιχα). Έτσι, ο *τόπος* σχετίζεται (χωρίς να ταυτίζεται) με τα όρια του σώματος που βρίσκεται εκεί, με τη *μορφή* ουσιαστικά του σώματος. Κατ' αντιστοιχία λοιπόν με την πλατωνική *Χώρα*, που δεν ταυτίζεται *ούτε* με τον Κόσμο των Ιδεών *ούτε* με τον Φυσικό Κόσμο, ο αριστοτελικός *τόπος* δεν είναι *ούτε* η μορφή των σωμάτων *ούτε* η ύλη τους.[15]

---

[12] «ὅμως τὸν τόπον καὶ τὴν χώραν τὸ αὐτὸ ἀπεφήνατο» Αριστοτέλης, *Φυσικά*, 209b15. Η διαφορά είναι ότι ο *Τόπος* για τον Αριστοτέλη δεν μπορεί να έχει αυτόνομη ύπαρξη.

[13] Για να εξηγήσει ορισμένες περιπτώσεις εξαναγκασμένης κίνησης, ο Αριστοτέλης είχε αναπτύξει τη θεωρία της *αντιπερίστασης*, σύμφωνα με την οποία η τάση ενός σώματος να κινείται ακόμα και όταν χάνει την επαφή του με το κινούν, οφείλεται σε εξωτερικά αίτια (π.χ. τον περιβάλλοντα αέρα).

[14] Ίσως με μια δόση ιστορικού αναχρονισμού, ο Jammer (2001, σ. 24) και άλλοι μελετητές, ερμηνεύουν τον αριστοτελικό χώρο όπως περιγράφεται στο χωρίο αυτό, ως «πεδίο δυνάμεων».

 Στη συνέχεια ο Αριστοτέλης υπεισέρχεται στο θέμα της ουσίας του χώρου. Η συζήτηση θα οδηγήσει σε έναν πρώτο ορισμό σύμφωνα με τον οποίο *τόπος* θεωρείται «τὸ πέρας τοῦ περιέχοντος σώματος» (*Φυσικά*, 212a 6), η κοίλη δηλαδή επιφάνεια του περιέχοντος σώματος και όχι το ίδιο το περιέχον σώμα. Ο Αριστοτέλης όμως δεν σταματά εκεί. Παραδεχόμενος ότι το να κατανοήσει κανείς την έννοια του *τόπου* είναι ένα δύσκολο ζήτημα, συγκρίνει τον *τόπο* με ένα αγγείο. Και όπως το αγγείο είναι *τόπος* φορητός, έτσι και ο *τόπος* είναι αγγείο αμετακίνητο (*Φυσικά*, 212a17). Κατά μία άποψη, στη φράση αυτή, απομονωμένη από το περικείμενό της, επικεντρώθηκε η πρόσληψη της αριστοτελικής θεωρίας του χώρου, ειδικά στον Μεσαίωνα (Μπετσάκος 2008).

 Μέσω της παραπάνω αντιπαραβολής ο Αριστοτέλης καταλήγει στον (τελικό) ορισμό του *τόπου*: «τὸ τοῦ περιέχοντος πέρας ἀκίνητον πρῶτον, τοῦτ' ἔστιν ὁ τόπος», θεωρώντας δηλαδή ότι ο *τόπος* είναι το πρώτο/ άμεσο *ακίνητο* όριο του περιέχοντος σώματος. Η θεώρηση αυτή είναι σύμφωνη με την ευρύτερη αριστοτελική φιλοσοφία, και επιπλέον μας δίνει τη δυνατότητα να δούμε από μια σύγχρονη οπτική γωνία τον χαρακτήρα του αριστοτελικού *τόπου*, ενδεχομένως και ως ακίνητο σύστημα συντεταγμένων (Γιανναράς 2001, σ. 310) απαραίτητο για τη μελέτη της κίνησης.

 Με βάση τα προηγούμενα, θα συμπληρώσουμε αυτή την παράγραφο παραθέτοντας την άποψη του Max Jammer σύμφωνα με την οποία τα αρχαία ελληνικά μαθηματικά δεν συνέβαλαν στη γεωμετρικοποίηση του χώρου, τουλάχιστον όσον αφορά στη χρήση χωρικών συντεταγμένων:

> «Ο χώρος, όπως υιοθετήθηκε από τη μηχανική ή την αστρονομία, δεν γεωμετρικοποιήθηκε ποτέ στην αρχαία ελληνική επιστήμη. Γιατί πώς θα ήταν δυνατόν ο ευκλείδειος χώρος, με την ομογένεια και την απειρία των γραμμών και των επιπέδων του, να ταιριάξει στο πεπερασμένο και ανισότροπο αριστοτελικό σύμπαν;».

---

[15] «οὐκ ἂν εἴη οὔτε ἡ ὕλη οὔτε τὸ εἶδος ὁ τόπος», *Φυσικά*, 210b29.

(Jammer 2001, σ. 34)

Είναι κάπως δύσκολο να γίνει άκριτα αποδεκτή η παραπάνω θέση, δεδομένου ότι το βασικό της επιχείρημα είναι η υποτιθέμενη ασυμβατότητα μεταξύ του *άπειρου* και *ομογενούς* ευκλείδειου χώρου από τη μια και του *περατού* και *ανισότροπου* αριστοτελικού σύμπαντος από την άλλη. Με δεδομένο το ευρύτερο φιλοσοφικό πλαίσιο της εποχής, τα έργα του Ευκλείδη και του Αριστοτέλη ίσως να μην είναι κατ' ανάγκη αλληλοαναιρούμενα, αλλά ενδεχομένως απόρροια διαφορετικών τρόπων θέασης της πραγματικότητας. Πράγματι, οι "λανθασμένες" απόψεις του Αριστοτέλη για τη λειτουργία του κόσμου, με τον τρόπο τους είναι "σωστές", υπό την έννοια ότι ο Σταγειρίτης έβλεπε τα πράγματα όχι τόσο διαφορετικά απ' ό,τι φαίνονται *από πρώτη ματιά*. Ο κόσμος του Αριστοτέλη διέπεται από *τάξη* η οποία δεν είναι άλλη από τον *τόπο* του εκάστοτε σώματος. Από τη στιγμή λοιπόν που αποδεχόμαστε την έννοια του *τόπου* ως τη φυσική θέση του σώματος, είναι λογικό να θεωρήσουμε ότι κάθε σώμα κινείται μέσα στον χώρο αναζητώντας τη θέση η οποία ανήκει από τη φύση του.[16] Υπ' αυτή την έννοια ο *τόπος* μπορεί να μην εναρμονίζει πλήρως, τείνει όμως να γεφυρώσει τα διεστώτα μεταξύ "ευκλείδειου χώρου" και "αριστοτελικού σύμπαντος".

Τέλος, υπάρχει σαφής ομοιότητα μεταξύ της περιγραφής των αποδεικτικών επιστημών στον Αριστοτέλη (*Αναλυτικά Ύστερα*, 72a 14-24 και 76a31 – 77a4) και του αξιωματικού συστήματος των *Στοιχείων* του Ευκλείδη. Η ομοιότητα αυτή οφείλεται κατά μείζονα λόγο στο γεγονός ότι οι κύριες έννοιες των *Στοιχείων* σφυρηλατήθηκαν σε μεγάλο βαθμό στην Ακαδημία του Πλάτωνα, όπου ο Αριστοτέλης έπαιξε πρωταγωνιστικό ρόλο. Για παράδειγμα, ο ορισμός της ευθείας στα *Στοιχεία* του Ευκλείδη, στον Πλάτωνα (*Παρμενίδης*, 137e3), τον Αριστοτέλη (*Τοπικά*, Ζ, 148b28) και τον Πρόκλο (*Σχόλια εις Ευκλείδην*, 110) είναι ουσιαστικά ο ίδιος και ταυτόχρονα μπορεί να αποτελέσει παράδειγμα πρωτογενούς γεωμετρικοποίησης, όπως στην περίπτωση κατά την οποία ο Ήλιος, η Σελήνη και η Γη βρίσκονται στην ίδια ευθεία (Σχήμα 3).

---

[16] Για το θέμα αυτό πολύ διαφωτιστική είναι η άποψη του T. Kuhn στο: Regis, E. (1995). *Ποιος πήρε την καρέκλα του Αϊνστάιν*. Μετ. Θ. Ηλιάδης. Αθήνα: Τροχαλία, κεφ. 9.

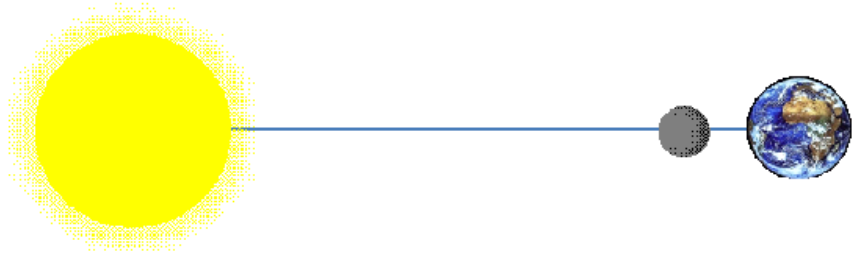

Σχήμα 3

### 3. Η επιρροή του Πλάτωνα και η αποθέωση του Αριστοτέλη

Από το 300 π.Χ. περίπου και για αρκετούς αιώνες, κέντρο της επιστημονικής δραστηριότητας υπήρξε αναμφίβολα η Αλεξάνδρεια. Ο Ευκλείδης, ο Απολλώνιος, ο Ίππαρχος, ο Ήρωνας, ο Πτολεμαίος, ο Διόφαντος, ο Πάππος και η Υπατία είναι μερικοί, μόνο, μαθηματικοί και αστρονόμοι που σχετίζονται με την Αλεξάνδρεια κατά την περίοδο αυτή. Η επιστημονική έρευνα, βέβαια, ακολούθησε (όχι πάντοτε γραμμικά) την κατεύθυνση που είχε υποδείξει η πλατωνική φιλοσοφία, σε συνδυασμό με την κρατούσα αριστοτελική φυσική. Έτσι, η ηλιοκεντρική θεωρία του Αρίσταρχου απορρίφθηκε γρήγορα από τον Απολλώνιο και τον Ίππαρχο οι οποίοι διατήρησαν το γεωκεντρικό δόγμα,[17] προετοιμάζοντας ταυτόχρονα το έδαφος για την επικράτηση του πτολεμαϊκού μοντέλου.

Στη συνέχεια οι νεοπλατωνικοί φιλόσοφοι (Πλωτίνος, Ιάμβλιχος, Πρόκλος κ.ά.) βασίστηκαν κυρίως στον *Τίμαιο* προκειμένου να εκφράσουν σε φιλοσοφικά και θρησκευτικά κείμενα τα κοσμολογικά τους μοντέλα ή τις ερμηνευτικές τους προσεγγίσεις. Η διδασκαλία τους είχε έναν χαρακτήρα "αποκάλυψης", επομένως ήταν πολύ φυσικό η γνήσια μαθηματική έρευνα και η παρατήρηση της φύσης να περάσουν σε δεύτερη μοίρα.

Ιδιαίτερη περίπτωση φιλοσόφου συνδεόμενου με τη νεοπλατωνική σχολή της Αλεξάνδρειας το πρώτο μισό του 6ου αιώνα, αποτελεί ο Ιωάννης Φιλόπονος, ο οποίος θεωρείται ο πρώτος εκπρόσωπος του χριστιανικού αριστοτελισμού. Με την επανερμηνεία φυσικών όρων και τη συγγραφή καίριων *σχολίων*,[18] άσκησε δημιουργική κριτική στο έργο του Αριστοτέλη προλειαίνοντας μέρος του δρόμου που οδήγησε στις νέες προσεγγίσεις των φυσικών επιστημών (Θεοδοσίου & Δανέζης 2010). Εστιάζοντας στις απόψεις του για τον χώρο και την κίνηση, θα λέγαμε πως ο Φιλόπονος πίστευε στην ύπαρξη ενός χώρου απόλυτου, γεμάτου με σώματα (Sambursky 1987, σ. 6), δεν απέκλειε την ύπαρξη του κενού, ενώ με νοητικά πειράματα αντέκρουσε τη θεωρία της αντιπερίστασης. Ερμήνευσε την κίνηση ενός σώματος που δεν βρισκόταν σε επαφή με το κινούν, με την εισαγωγή μιας ωθητικής δύναμης την οποία αργότερα ο Jean Buridan ονόμασε *impetus*. Κατά τη διάρκεια του ύστερου Μεσαίωνα η θεωρία αυτή γνώρισε μεγάλη απήχηση στη Δυτική Ευρώπη (Αραμπατζής κ.ά. 1999).

---

[17] Για το θέμα αυτό βλ. Καρτσωνάκης, Μ. (1992). "Αιτίες για τη μη αποδοχή της ηλιοκεντρικής θεώρησης του Αρίσταρχου: Όψεις μιας αστρονομικής διαμάχης", στο Δ. Α. Αναπολιτάνος – Β. Καρασμάνης (επιμ.), *Κείμενα Ιστορίας και Φιλοσοφίας των Αρχαίων Ελληνικών Μαθηματικών*, σσ. 93-110, Αθήνα: Τροχαλία.

Έτσι, χωρίς να υπάρχει σαφής διαχωριστική γραμμή, φτάνουμε στο «προαύλιο του βυζαντινού κόσμου» το οποίο χαρακτηρίζεται από φιλοσοφικές διαμάχες και θρησκευτικές αντιπαλότητες, όπως κάθε περίοδος κρίσης. Αυτή η πνευματική ζύμωση, όμως, συνέβαλε στο να παραμείνει ζωντανή η επιστημονική παράδοση, και η ουσία του ελληνικού πνεύματος να περάσει στον χώρο του ανατολικού ρωμαϊκού κράτους. Ωστόσο από όλες τις φιλοσοφικές σχολές και τάσεις του κλασικού ελληνισμού, μόνο ο πλατωνισμός και ο αριστοτελισμός ξεπέρασαν την ύστερη αρχαιότητα και άσκησαν μια μόνιμη, λίγο-πολύ, επίδραση στην πνευματική ιστορία του Μεσαίωνα (Hunger 2004, τ. Α΄, σ. 50). Για παράδειγμα, η μελέτη του *Τίμαιου* και των *Φυσικών* και η επιχειρούμενη σύνθεσή τους σε μία θρησκευτική (και παιδαγωγικά προσιτή) βάση, οδήγησε μεγάλες πνευματικές μορφές του Βυζαντίου όπως τον Μιχαήλ Ψελλό,[19] τον Νικηφόρο Βλεμμύδη,[20] τον Θεόδωρο Λάσκαρι[21] και τον Πλήθωνα Γεμιστό,[22] να συγγράψουν σχό-λια και πρωτότυπα έργα τα οποία αποτέλεσαν διδακτικά εγχειρίδια τόσο για την Ανατολή όσο και για τη Δύση.

---

[18] Η εκτεταμένη γραμματεία των αριστοτελικών σχολιαστών (Αλέξανδρος Αφροδισιέας, Θεμίστιος, Σιμπλίκιος, Ιωάννης Φιλόπονος) δείχνει ότι στην ύστερη αρχαιότητα ο Αριστοτέλης θεωρείται σταθμός στην ιστορία της φιλοσοφίας (Κάλφας 2014, σ. 144).

[19] Βλ. τα *Σχόλια εις τα Φυσικά του Αριστοτέλους* στο Psellos, M. (2008). *Kommentar zur Physik des Aristoteles*. Editio princeps. Einleitung, Text, Indices von L. G. Benakis. Αθήναι: Ακαδημία Αθηνών.

[20] Νικηφόρου του Βλεμμύδου, (1863). "Τὰ εὑρισκόμενα πάντα", στο J. P. Migne (επιμ.), *Patrologia Graeca*, 142.

[21] Βλ. Pappadopulos, J. (1908). *Théodore II Laskaris, Empereur de Nicée*, Paris: Librairie Alphonse Picard et fils και Ivánka, E. (1972). "Mathematische Symbolik in den beiden Schriften des Kaisers Theodoros II Laskaris: Δήλωσις φυσική und Περὶ φυσικῆς κοινωνίας", *Byz. Forsch*, 4, pp. 138-141.

[22] Γεωργίου Γεμιστού του και Πλήθωνος, (1866). "Περὶ ὧν Ἀριστοτέλης πρὸς Πλάτωνα διαφέρεται", στο J. P. Migne (επιμ.), *Patrologia Graeca*, 160, pp. 889-934.

Από τον 11ο αιώνα, στο πλαίσιο του σχολαστικισμού, τα έργα του Πλάτωνα και του Αριστοτέλη άρχισαν να μελετώνται σταδιακά και στη Δυτική Ευρώπη, αφού έως τότε τα κείμενα των ελλήνων φιλοσόφων είτε δεν ήταν διαθέσιμα, είτε ήταν γνωστά μόνο με αποσπασματικό τρόπο. Βαθμιαία, όμως, όλο και περισσότερες λατινικές μεταφράσεις έκαναν την εμφάνισή τους, ιδίως τον 12ο αιώνα, ο οποίος ονομάστηκε «εποχή της μετάφρασης».[23] Έτσι, η δυτική μεσαιωνική επιστήμη αναπτύχθηκε βασισμένη στην αριστοτελική κοσμοθεώρηση, ενίοτε όμως και εν μέσω αντι-αριστοτελικών κριτικών, είτε λόγω έλλειψης ικανοποιητικών εξηγήσεων, είτε στη βάση νέων θεολογικών ερμηνειών ή αστρονομικών δεδομένων, σε κάθε περίπτωση πάντως, με σημείο αναφοράς τον Σταγειρίτη φιλόσοφο. Η συνεχής αυτή κριτική οδήγησε τελικά, τον 17ο αιώνα, στη διατύπωση επιστημονικών θεωριών, οι οποίες προέκριναν ορθολογικά εναλλακτικά μοντέλα για τον χώρο και την εξήγηση των φυσικών φαινομένων.

Ταυτόχρονα, η επιρροή του Πλάτωνα αν και εξασθενημένη, εξακολουθούσε να υπάρχει στον ελληνικό και τον δυτικό κόσμο. Εξάλλου, η παράδοση του βυζαντινού πλατωνισμού ποτέ δεν διακόπηκε εντελώς από την εποχή που έκλεισε η Ακαδημία των Αθηνών έως την εποχή που μια νέα Ακαδημία ιδρύθηκε το 1459 στη Φλωρεντία (Klibansky 1939) από τον Πλήθωνα Γεμιστό με την υποστήριξη του Κοσμά των Μεδίκων (Cosimo de' Medici). Στα σπουδαιότερα κέντρα μάθησης της Δύσης ο *Τίμαιος*, ο πιο γνωστός πλατωνικός διάλογος στον Μεσαίωνα (Hunger 2004, σ. 52), μελετήθηκε εντατικά και χρησιμοποιήθηκε για να εξηγηθεί η δομή του σύμπαντος (Grant 1994, σ. 24). Ας μην ξεχνάμε άλλωστε ότι τα πέντε κανονικά στερεά είχαν δημιουργήσει ένα κλίμα μυστικισμού καθ' όλη τη διάρκεια του Μεσαίωνα, μέχρι και το 1596 όπου εμφανίστηκαν για τελευταία φορά στο *Mysterium Cosmographicum* του Kepler, και υποτίθεται ότι κα-θόριζαν τις κυκλικές τροχιές των πλανητών, ως εγγραφόμενα σε αντίστοιχες ομόκεντρες σφαίρες.

---

[23] Στη μεταφορά της αριστοτελικής φυσικής φιλοσοφίας στη Δυτική Ευρώπη, φαίνεται ότι συνέβαλε η πρόσληψη του Αριστοτέλη από τους Άραβες, των οποίων η φιλοσοφία και επιστήμη είχε από την αρχή δείξει προτίμηση στον Σταγειρίτη. Το πρόγραμμα μεταφράσεων που είχαν αναπτύξει, σκόπευε να κάνει γνωστά τα έργα του «φιλοσόφου», όπως οι ίδιοι αποκαλούσαν τον Αριστοτέλη, στη Δύση. Ωστόσο, νεώτεροι ιστορικοί υποστηρίζουν ότι η αρχαία ελληνική γραμματεία, και εν πολλοίς και ο Αριστοτέλης, μεταφέρθηκαν στη Δύση από πολύ νωρίς και κυρίως χάρη στους Βυζαντινούς και εν συνεχεία από τους Άραβες και όχι το αντίθετο. Βλ. σχετ. Lemerle, P. (2007). *Ο πρώτος Βυζαντινός ουμανισμός*. Μετ. Μ. Νυσταζοπούλου-Πελεκίδου. Αθήνα: Μ.Ι.Ε.Τ. και Gouguenheim, S. (2009). *Ο Αριστοτέλης στο Μον-Σαιν Μισέλ*. Μετ. Φ. Γαϊδατζή – Φ. Μπούμπουλη. Αθήνα: Ολκός.

### 4. Το σύμπαν ξαναγίνεται άπειρο

Σχεδόν 2.000 χρόνια μετά τους Ίωνες φιλοσόφους, η ιδέα ενός απέραντου χώρου έρχεται και πάλι στο προσκήνιο από τον Nicola Cusano (ή Niccolò da Cusa) με το έργο του *De docta ignorantia* (*Η πεφωτισμένη άγνοια*) στα 1440. Η πραγματεία αυτή, ακροβατώντας μεταξύ Θεολογίας και Αστρονομίας, προβάλλει ως ένα είδος *νοητικού πειράματος* προορισμένου να συμβάλει στην υπέρβαση των περιορισμών που θέτουν η ορθολογική σκέψη και η ατελής φύση της ανθρώπινης γνώσης. Εκεί μπορεί κανείς να διακρίνει στοιχεία που αφορούν στη σχετικότητα των κινήσεων (π.χ. αργό – γρήγορο) και του χώρου (π.χ. προσανατολισμός), καθώς και στη δυνατότητα της συνεχούς μεταβολής ενός σχήματος ως προς τη μορφή και το μέγεθός του. Τα σημαντικότερα, ίσως, στοιχεία της κοσμολογίας του Cusano είναι η απόρριψη της ιεραρχικής δομής του σύμπαντος και της κεντρικής θέσης που υποτίθεται ότι σ' αυτό κατέχει η Γη. Δυστυχώς, τα συμπεράσματα αυτά, βασισμένα κυρίως στη διαίσθηση, υπονομεύτηκαν από τις επιστημονικές του απόψεις, οι οποίες συνιστούσαν μάλλον οπισθοδρόμηση (Koyré 1989, σ. 27).

Ο πρώτος, όμως, που πίστεψε και διακήρυξε με ενθουσιασμό την ύπαρξη ενός άπειρου σύμπαντος, ενός ομογενούς χώρου δίχως κέντρο, ήταν ο Giordano Bruno. Επηρεασμένος από την πραγματεία του Λουκρήτιου *De rerum natura* όπου κυριαρχούν επιχειρήματα για την απειρία του σύμπαντος, είναι πια σε θέση να επιτεθεί στον πεπερασμένο αριστοτελικό κόσμο. Έτσι, στο έργο του *De l' infinito, universo et Mondi*, ο Bruno "συνδιαλέγεται" με τον Αριστοτέλη, προκειμένου να αποδυναμώσει τα επιχειρήματά του –κυρίως αυτά που αφορούν στην *τάξη* των σωμάτων– υιοθετώντας μια πρωτοπόρα για την εποχή του συλλογιστική, όπως για παράδειγμα τις αναφορές σε μαθηματικό και φυσικό χώρο, την ύπαρξη και άλλων (άπειρων) κόσμων ή την απουσία εξώτατης σφαίρας και πρώτου κινούντος.

Τα πράγματα, πλέον, παίρνουν τον δρόμο τους. Οι ουράνιες σφαίρες σπάνε, ο κόσμος παύει να είναι κλειστός. Ο χώρος απελευθερώνεται από το εσωτερικό ενός πεπερασμένου σύμπαντος· γίνεται άπειρος και ομογενής. Οι περιορισμοί που έθεταν οι ανθρώπινες αισθήσεις, οι κοινωνικές προκαταλήψεις και οι θεολογικές προκείμενες καταργούνται. Η Γη παύει να είναι το κέντρο του κόσμου. Ο χώρος γεωμετρικοποιείται, γίνεται Ευκλείδειος. Για να αντικατασταθεί, όμως, πλήρως το Πτολεμαϊκό μοντέλο από ένα σύγ-χρονο ηλιοκεντρικό, θα χρειαστούν οι αστρονομικές παρατηρήσεις του T. Brahe και η διάνοια του J. Kepler.

### 5. Μια (όχι και τόσο) γνωστή διαμάχη για τη φύση του χώρου

Στην κλασική πραγματεία του, *Μαθηματικές Αρχές της Φυσικής Φιλοσοφίας* (*Philosophiae Naturalis Principia Mathematica*), τη βάση της κλασικής Φυσικής, ο Νεύτωνας με την εισαγωγή των κατάλληλων μαθηματικών εργαλείων ουσιαστικά διαχωρίζει τη φυσική έρευνα από τη μεταφυσική και τη Θεολογία. Αν και το έργο επιδέχεται διαφορετικές επιστημολογικές και φιλοσοφικές ερμηνείες, εντούτοις όλα συνηγορούν στο ότι για τον Νεύτωνα ο *χώρος* και ο *χρόνος* αποτελούν απροσδιόριστες και ενδεχομένως έννοιες με μεταφυσική χροιά, σε τέτοιο βαθμό που η φυσική του φιλοσοφία να μην μπορεί να διαχωριστεί από αυτές. Όπως λέει ο ίδιος:

«Δεν ορίζω τον χρόνο, τον χώρο, τον τόπο και την κίνηση, ως πολύ γνωστά σε όλους. Πρέπει μόνο να παρατηρήσω ότι το ευρύ κοινό αντιλαμβάνεται αυτές τις ποσότητες μόνο αναφορικά με τα αισθητά αντικείμενα. Έτσι, προκύπτουν ορισμένες προκαταλήψεις, για την άρση των οποίων θα ήταν χρήσιμο να διακρίνουμε τις έννοιες αυτές σε απόλυτες και σχετικές, αληθείς και φαινόμενες, μαθηματικές και κοινές».
(Newton 1962, σ. 6, δική μας απόδοση)

Σημαντικό ρόλο στη θεώρηση αυτή του Νεύτωνα έχει σίγουρα παίξει το θρησκευτικό στοιχείο, το οποίο διέπει κυρίως το ύστερο έργο του, καθώς και το κοινωνικό πλαίσιο της εποχής του με τις αλλεπάλληλες πολιτικές διαμάχες και τις επιστημονικές καινοτομίες.[24] Όλα αυτά συνετέλεσαν στο να ταυτιστεί ο απόλυτος χώρος[25] με τον Θεό ή με την απανταχού παρουσία Του. Έτσι, ο κόσμος του Νεύτωνα είναι μια μηχανή που έχει κατασκευαστεί και ρυθμιστεί από τον Θεό αλλά από καιρού εις καιρόν χρειάζεται κάποιες διορθωτικές παρεμβάσεις υπό μορφή ανανέωσης του ενεργειακού αποθέματος για να διατηρείται η λειτουργία της, η οποία από ορισμένες ανεξήγητες πλανητικές κινήσεις κινδυνεύει να διαταραχθεί.

Η παραπάνω θέση είναι μια από τις κύριες εστίες αντεγκλήσεων με τον Leibniz: ο Γερμανός φιλόσοφος ήταν οπτιμιστής, είχε δηλαδή την πεποίθηση ότι ο κόσμος μας, ως έκφραση της *Γνώσης*, της *Βούλησης* και της

---

[24] Για το ιστορικό πλαίσιο βλ. Hill, C. (2011). *Ο Αιώνας της Επανάστασης 1603-1714*. Μετ. Ελένη Αστερίου. Αθήνα: Οδυσσέας.
[25] Ο Νευτώνειος χώρος είναι ένας ουσιαστικά ευκλείδειος χώρος με τρεις απόλυτες διαστάσεις, άπειρος και ακίνητος.

*Ισχύος* του Θεού, είναι ο καλύτερος δυνατός κόσμος. Έτσι, ενώ ο Νεύτωνας θεωρεί ότι ο Θεός επέλεξε (ίσως τυχαία) να δημιουργήσει *αυτόν* τον κόσμο ανάμεσα σε πολλούς άλλους, ο Leibniz "επιβάλλει" στον Θεό έναν ορισμένο τρόπο δράσης, κάτι δηλαδή που, σύμφωνα με τους επικριτές του, οδηγεί στην αναγκαιότητα. Άλλωστε η έριδα μεταξύ Νεύτωνα και Leibniz ήταν κατά βάση μια σύγκρουση φιλοσοφικών κοσμοθεωρήσεων για τη φύση του Θεού, της ύλης και της δύναμης (Iltis 1973).

Το φιλοσοφικό σύστημα του Leibniz (όπως και αυτό των Descartes και Spinoza) ανήκει στην ορθολογική φιλοσοφική παράδοση του 17$^{ου}$ αιώνα. Ο Leibniz φαίνεται να δέχεται την ύπαρξη τριών επιπέδων, τριών Κόσμων θα λέγαμε, στους οποίους διαιρείται το οντολογικό του σχήμα.[26] Στο πρώτο επίπεδο, το επίπεδο του *Πραγματικού*, υπάρχουν οι *Μονάδες*, δηλαδή απλές ουσίες, οι οποίες είναι αδιάστατες, αδιαίρετες και άφθαρτες. Τα *όντα* αυτά παίζουν τον ρόλο των δομικών στοιχείων του σύμπαντος, δεν καταλαμβάνουν όμως κάποιον *χώρο* και δεν μεταβάλλονται με την πάροδο του *χρόνου*. Στο δεύτερο επίπεδο, αυτό των *Φαινομένων*, ανήκει η χωροχρονικότητα και τα αντικείμενά της. Εδώ φαίνεται ότι ο *χώρος* και κατ' επέκταση ο *χρόνος* είναι για τον Leibniz *συστήματα σχέσεων μεταξύ υλικών αντικειμένων*, απαλλαγμένα από μεταφυσικά στοιχεία. Ο χώρος είναι ο τρόπος συνύπαρξης των σωμάτων, ενώ ο χρόνος είναι η σειρά διαδοχής των συμβάντων. Είναι «το σύνολο των χωροχρονικών σχέσεων που γίνονται αντιληπτές, ως τέτοιες, στο πλαίσιο της εσωτερικής προοπτικής θέασης των αναπαραστασιακών καταστάσεων των Μονάδων».[27] Στο τρίτο επίπεδο, τέλος, το επίπεδο του *Ιδεώδους*, υπάρχουν αφηρημένες ιδέες και καθολικές αλήθειες, οι οποίες προϋπάρχουν στον νου των λογικών όντων και με την εμπειρία απλώς αναδύονται, κάτι το οποίο μας παραπέμπει ευθέως στον Πλάτωνα.

Κεντρικό ρόλο στο φιλοσοφικό οικοδόμημα του Leibniz κατέχουν δύο Αρχές, οι οποίες ανάγονται στον Αριστοτέλη: η *Αρχή της Αντιφάσεως*[28] και η *Αρχή του Αποχρώντος Λόγου*. Οι Αρχές αυτές οδηγούν, κατά κάποιο τρόπο, τον Leibniz να "ορίσει" τον *χώρο* και τον *χρόνο* με τον τρόπο που είδαμε προηγουμένως, ενώ αντίθετα δίνουν την ευκαιρία στους φιλοσόφους των άλλων ρευμάτων να του προσάψουν ότι ο Αποχρών Λόγος μπορεί να είναι απλά η θέληση του Θεού και όχι κάποια φυσική αιτία. Όλα αυτά λοιπόν, ωθούν σχεδόν αναπόφευκτα τον Leibniz να συνδέσει τον *χώρο* με τον *χρόνο*, μια θέση δηλαδή η οποία βρίσκεται πιο κοντά στις σύγχρονες σχετικιστικές αντιλήψεις, σε αντίθεση με τον Νεύτωνα που θεωρούσε τον χώρο και τον χρόνο δύο ανεξάρτητες και απόλυτες έννοιες.

---

[26] Βλ. την εισαγωγή του επιμελητή στο Leibniz, G. W. (2006). *Η Μοναδολογία*. Μετάφραση Στ. Λαζαρίδης, εισαγωγή-επιμέλεια Δ. Αναπολιτάνος. Αθήνα: Εκκρεμές.
[27] Leibniz, *ό.π.* σ. 24.
[28] Είναι πολύ ενδιαφέρουσα η διατύπωση, από τον ίδιο τον Leibniz, της Αρχής της Αντιφάσεως, που εμπεριέχει έναν "χωρικό" χαρακτηρισμό του ψευδούς: «… nous jugeons *faux* ce qui en enveloppe» (… κρίνουμε *ψευδές* αυτό που αναδιπλώνεται στον εαυτό του). Leibniz, *ό.π.* σσ. 50-51, δική μας απόδοση.

### 6. Από τον χώρο και τον χρόνο στον χωροχρόνο

Η έννοια του απόλυτου χώρου του Νεύτωνα, τουλάχιστον όσον αφορά στα Μαθηματικά, διατηρήθηκε μέχρι την ανακάλυψη των μη-Ευκλείδειων γεωμετριών, ενώ η έννοια του απόλυτου χρόνου συντηρήθηκε για λίγα ακόμη χρόνια. Στη Φυσική, μάλιστα, προκειμένου να "σωθούν τα φαινόμενα" –να διατηρηθούν αναλλοίωτες οι εξισώσεις Maxwell κάτω από τους μετασχηματισμούς Lorentz–, προτάθηκε ο απόλυτος χώρος να ταυτιστεί με τον αιθέρα, θεωρούμενο ως φορέα των ηλεκτρομαγνητικών κυμάτων και ταυτόχρονα ως προνομιακό σύστημα αναφοράς. Η πρόταση αυτή εγκαταλείφθηκε οριστικά με τη διατύπωση της Ειδικής Θεωρίας της Σχετικότητας, στην οποία παρά ταύτα ο χώρος παρέμενε μια αινιγματική έννοια:

«*Πρώτα απ' όλα ας αποφύγουμε την ασαφή λέξη "χώρος" με την οποία, πρέπει έντιμα να ομολογήσουμε, δεν μπορούμε απολύτως τίποτα να κατανοήσουμε*».
(Einstein 2015, σ. 18, σε δική μας απόδοση)

Αυτές οι επαναστατικές αλλαγές στην επιστημονική και φιλοσοφική θεώρηση που επέφερε η διατύπωση της Θεωρίας της Σχετικότητας από τον Einstein, αναλύθηκαν διεξοδικά από τον Weyl (1952) μαζί με τη σε βάθος μελέτη των εννοιών του χώρου, του χρόνου και της ύλης.

Με το τετραδιάστατο χωροχρονικό συνεχές ο "χώρος" πλέον μετατρέπεται σε ένα αφηρημένο μαθηματικό αντικείμενο, πληρούται με σημειογεγονότα, αποκτά άλλη υπόσταση και δική του αναπαράσταση με τη βοήθεια του *κώνου φωτός* ή του *διαγράμματος Minkowski*. Όπως σημειώνει ο δημιουργός του:

«*Οι θεωρήσεις για τον χώρο και τον χρόνο τις οποίες θέλω να θέσω υπόψη σας, ξεπήδησαν από το πεδίο της πειραματικής φυσικής, και από εκεί εξάπλωσαν τη δύναμή τους. Είναι ριζοσπαστικές. Εφεξής ο χώρος και ο χρόνος ξεχωριστά ο ένας από τον άλλο, είναι καταδικασμένοι να καταντήσουν απλές σκιές, και μονάχα ένα είδος ένωσης και των δύο θα διατηρήσει μια αυθύπαρκτη υπόσταση*».
(Minkowski 1952, σ. 75, δική μας απόδοση και έμφαση)

Για να ολοκληρωθεί η παραπάνω ιδέα απομένει η μαθηματική περιγραφή του κώνου φωτός. Αυτή βασίζεται στην ταξινόμηση των γεγονότων ως προς ένα αυθαίρετο γεγονός P σε τρεις κατηγορίες, ανάλογα με το αν το τετράγωνο του διανύσματος μετατόπισής τους από το P είναι θετικό, αρνητικό ή μηδέν. Χωρίς να παραθέτουμε τεχνικές λεπτομέρειες, σημειώνουμε απλά πως κατ' αυτό τον τρόπο ορίζεται ένας γεωμετρικός τόπος, που είναι ουσιαστικά ένας κώνος φωτός, ανεξάρτητα από το πλήθος των διαστάσεων.

### 7. Αντί επιλόγου

Από τις αρχές της δημιουργίας της ως επιστημονικού κλάδου, η Γεωμετρία είναι θεμελιωμένη επάνω σε αξιώματα που αναφέρονται σε "σημεία", "ευθείες", "επίπεδα" κ.λπ. Το ερώτημα που ανακύπτει εδώ έχει να κάνει με τη *σχέση* αυτών των οντοτήτων ως προς την εμπειρική γνώση που προκύπτει από την παρατήρηση. Σύμφωνα με τον Πουανκαρέ:

> «[…] δεν πειραματιζόμαστε επί ευθειών ή κύκλων ιδανικών. Μπορούμε να πραγματοποιήσουμε πειράματα μόνο επί υλικών αντικειμένων».

(Poincaré 1952, σ. 49, δική μας απόδοση)

Η άποψη αυτή δείχνει τη διάσταση μεταξύ μιας αξιωματικά θεμελιωμένης Γεωμετρίας, με αξιώματα δοσμένα a priori, και μιας εξέλιξης κατευθυνόμενης από το πείραμα η οποία επιβάλλει περιορισμούς και επιτρέπει αναθεωρητικές έρευνες «επί των υποθέσεων που κείνται στα θεμέλια της Γεωμετρίας», όπως υποδεικνύει και ο τίτλος της διατριβής του Ρίμαν (βλ. Riemann 2016).

Δεδομένου, λοιπόν, ότι μέχρι τους Gauss, Lobachevsky και Bolyai ο χώρος ήταν μια απόλυτη οντότητα και εθεωρείτο Ευκλείδειος, η ανακάλυψη των μη-Ευκλείδειων γεωμετριών είχε επίπτωση στην επιστημολογία του χώρου "επιβάλλοντας" την έννοια του *μαθηματικού χώρου*. Στη συνέχεια η Ειδική και κυρίως η Γενική Θεωρία της Σχετικότητας "επέβαλαν" την έννοια του *φυσικού χώρου*, ώστε η σύγχρονη Φυσική να βασίζει, πλέον, τη θεώρησή της για τον χώρο στη ρημάνεια έννοια της πολλαπλότητας.

Συνεπώς η έννοια της ριμάνειας πολλαπλότητας και εν συνεχεία η Γεωμετρία Minkowski είχαν ως αποτέλεσμα την πληρέστερη κατανόηση της υποθετικής φύσης των γεωμετρικών αξιωμάτων και την υπόδειξη της απουσίας μιας γεωμετρίας a priori κατάλληλης να περιγράψει χωρικές σχέσεις μεταξύ φυσικών σωμάτων. Επιπλέον, με την εισαγωγή του τετραδιάστατου χωροχρονικού συνεχούς, ο χώρος έπαψε να είναι απόλυτος (π.χ. ο χώρος του Νεύτωνα ή ο αιθέρας του Lorentz) και άνοιξε ο δρόμος για την υιοθέτηση της έννοιας του πεδίου και τη διατύπωση της Γενικής Θεωρίας της Σχετικότητας. Ο δρόμος αυτός οδηγεί στην αιχμή της επιστημονικής έρευνας σήμερα, στα *βαρυτικά κύματα* και τη διασύνδεσή τους (;) με το πλήθος των διαστάσεων του σύμπαντος.